%

\magnification=\magstep1
\input amstex
\documentstyle{amsppt}

\topmatter
\title
Semi-Cohen versus Cohen algebras 
\endtitle
\author
Jind\v rich Zapletal
\endauthor
\thanks
\endthanks
\affil
The Pennsylvania State University
\endaffil
\address
Department of Mathematics, 
The Pennsylvania State University, 
University Park,  PA 16802
\endaddress
\email
zapletal\@math.psu.edu
\endemail
\abstract
We show that there are semi-Cohen Boolean algebras which cannot be 
completely embedded into Cohen Boolean algebras. Using the ideas
 from this proof, we give a simpler argument for a theorem of S. 
Koppelberg and S. Shelah, stating that there are complete subalgebras 
of Cohen algebras which are not Cohen themselves.
\endabstract
\endtopmatter

\document
\subhead
{0. Introduction}
\endsubhead

Probably the best-known and most analyzed species of complete Boolean 
algebras are the Cohen algebras $\Bbb C(\kappa ),$ representing the 
regular open subsets of $\{ 0,1\} ^\kappa .$ Their internal structure 
is quite well understood from the algebraic as well as forcing point 
of view. Motivated by some previous work \cite {Ko}, S. Fuchino and T. 
Jech proposed the following notion, intended to capture all of the 
combinatorial content of ``being a $\Bbb C(\kappa )":$

\proclaim
{Definition 1} A complete Boolean algebra $B$ is {\it semi-Cohen} if 
$S=\{ A\in [B]^{\aleph _0}:A$ is a regular subalgebra of $B\}$ contains 
a closed unbounded set in $[B]^{\aleph _0}.$
\endproclaim

Thus semi-Cohen algebras behave locally much like Cohen algebras. 
This would suggest a close relationship, if not equality, between 
the two classes of algebras. Indeed, it turns out that every 
$\Bbb C(\kappa )$ is semi-Cohen and every semi-Cohen algebra 
of uniform density $\aleph _1$ is isomorphic to $\Bbb C(\aleph _1).$ We have

\proclaim
{Question 1} Is every semi-Cohen algebra of uniform density $\kappa$ 
isomorphic to $\Bbb C(\kappa )?$
\endproclaim

Now being semi-Cohen is a property hereditary to complete subalgebras, 
and it has been unclear how to obtain a similar result for Cohen algebras. 
Thus, the following two approximations to Question 1 are interesting, too:

\proclaim
{Question 2} \cite {J2} Does every semi-Cohen algebra completely embed 
into a Cohen algebra?
\endproclaim

\proclaim
{Question 3} \cite {J2, Ka, Ko, KS} Is every complete subalgebra of 
$\Bbb C(\kappa )$ of uniform density $\kappa$ isomorphic to $\Bbb C(\kappa )?$
\endproclaim

Recently, S. Koppelberg and S. Shelah \cite {KS} showed that Question 3, 
and {\it a posteriori} Question 1, has a negative answer. We prove that 
Question 2 has a negative answer as well and provide an easier argument
 for the theorem in \cite {KS}.In the first section we give an example 
of a semi-Cohen algebra of uniform density $\frak c^+$ which cannot be 
embedded into a Cohen algebra. In the second section we exhibit a 
complete subalgebra of $\Bbb C(\aleph _2 )$ of uniform density $\aleph _2$ 
which is not isomorphic to $\Bbb C(\aleph _2 ).$

Theorem 1 is due to the author. Theorem 2 is due to  S. Koppelberg and 
S. Shelah \cite {KS}; we present a much cleaner and simpler argument. 
The author would like to thank B. Balcar and T. Jech for patiently 
listening to the proofs.

Our notation follows the set-theoretic standard as set forth in \cite 
{J1}. In a forcing notion we write $p\geq q$ to mean that $q$ is 
more informative than $p.$  also $p\perp q$ means that $p$ and $q$ 
are incompatible. If $\langle P,\leq \rangle$ is a partially ordered set 
(a {\it poset}) and $R\subset P$ we say that $R$ is a {\it regular subposet} 
of $P$ if every maximal antichain in $\langle R, \leq \rangle$ is a 
maximal antichain in $\langle P, \leq \rangle ,$ or, equivalently, 
$\forall p_0\in P\ \exists p_1\in R\ \forall p_2\in R\ p_2\leq p_1$ 
implies that $p_2$ and $p_0$ are compatible. For $I\subset Ord$ we 
understand $\Bbb C(I)$ as $RO(C(I)),$ where $C(I)=\{ p:p$ is a 
function, $dom(p)\in  [I]^{<\aleph _0}$ and $rng(p)\subset 2 \}$ 
ordered by extension. $H_\theta$ is a collection of sets hereditarily 
of cardinality $<\theta .$ $TC(x)$ is the transitive closure of $x.$
 We define a strict ordering $\ll$ on ${^\omega \omega}$ by $f\ll g$ 
if (a) for all but finitely many $n\in \omega \ f(n)\leq g(n)$ and (b) 
for infinitely many $n\in \omega \ f(n)<g(n).$ 

\subhead
{1. The first example}
\endsubhead

We give a negative answer to Question 2.

\proclaim
{Theorem 1} There is a semi-Cohen algebra of uniform density $\frak c^+$ 
which cannot be embedded into a Cohen algebra.
\endproclaim

\demo
{Proof} As usual, we prefer to deal with posets rather than with 
complete Boolean algebras. In this spirit, we give

\proclaim
{Definition 2} A poset $\langle P,\leq \rangle$ is semi-Cohen if 
$S=\{ R\in [P]^{\aleph _0}: R$ is a regular subposet of $P\}$ 
contains a closed unbounded set in $[P]^{\aleph _0}.$
\endproclaim

We can easily justify our choice of terminology by

\proclaim
{Lemma 1} $P$ is semi-Cohen iff $RO(P)$ is semi-Cohen. \qed L1
\endproclaim

Let $I\subset Ord.$ We define a poset $P(I)$ as follows:

\proclaim
{Definition 3} $P(I)=\{ p:p$ is a function, $dom(p)\in [I]^{<\aleph _0},
\exists n_p \in \omega \ \forall \alpha \in dom(p)$ $p(\alpha )\in {^{n_p}
\omega }\}$ ordered by $p_0\geq p_1$ if $dom(p_0)\subset dom(p_1),$ 
$n_{p_0}\leq n_{p_1},$ $\forall \alpha \in dom(p_0)$ $p_0(\alpha )
\subset p_1(\alpha )$ and $\forall i\in n_{p_1}\setminus n_{p_0}\ 
\forall \beta <\alpha$ both in $dom(p_0)$ we have $p_1(\beta )(i)
\leq p_1(\alpha )(i).$
\endproclaim

The important properties of the $P(I)$'s are recorded in

\proclaim
{Lemma 2}\roster
\item If $J\subset I$ then $P(J)$ is a regular subposet of $P(I);$ 
consequently, $P(I)$ is semi-Cohen.
\item $P(I)$ has uniform density $|I|+\aleph _0.$
\item Let $F\subset P$ be generic. For $\alpha \in I$ we define 
$f_\alpha =\bigcup _{p\in F}p(\alpha )\in {^\omega \omega }.$ 
Then $\alpha <\beta \in I$ implies $f_\alpha \ll f_\beta .$
\item $\langle f_\alpha :\alpha \in I\rangle$ determines $F:$ 
namely, $G=\{ p\in P(I):\forall \alpha \in dom(p)\ p(\alpha )
\subset f_\alpha$ and for all $\alpha >\beta$ both in $dom(p)$ 
and $i\in \omega \setminus n_p$ $f_\beta (i)\leq f_\alpha (i)\} .$ 
\endroster
\endproclaim

\demo
{Proof} We start with (1). Let $J\subset I$ and $p_0\in P(I).$ 
Set $p_1=p_0\restriction J.$ We show that $\forall p_2\in P_j$ 
$p_2\leq p_1$ implies that $p_0,p_2$ are compatible, finishing 
the proof of regularity. To this aim, fix $p_2\in P(J), p_2\leq p_1.$ 
Without loss of generality we may assume that $n_{p_2}\geq n_{p_0}.$ 
(If this is not the case, just strengthen $p_2$ within $P(J)$ to make 
it hold.) We define $p_3\in P(I),$ a common lower bound of $p_0,p_2.$ 
We let $dom(p_3)=dom(p_0)\cup dom(p_2),$ $n_{p_3}=n_{p_2},$ for $\alpha 
\in dom(p_2)$ $p_3(\alpha )=p_2(\alpha )$ and for $\alpha \in dom(p_0)
\setminus J$ we set $p_3(\alpha )$ to be the function in ${^{n_{p_2}}
\omega }$ extending $p_0(\alpha )$ such that $\forall n\in n_{p_2}
\setminus n_{p_0}$ $p_3(\alpha )(n)=p_2(\beta )(n),$ if $dom(p_0)\cap 
J\cap \alpha \neq 0$ and $\beta =max(dom(p_0)\cap J\cap \alpha ),$ 
and $p_3(\alpha )(n)=0$ otherwise. By the definition of $P(I),$ we have
 $p_3\leq p_2,p_3\leq p_0$ and we are finished.

The immediate consequence of this is that $P(I)$ is semi-Cohen: 
$S=\{ P(J):J\in [I]^{\aleph _0}\}$ is a closed unbounded set in 
$[P]^{\aleph _0}$ consisting of regular subposets of $P,$ as we have just proven.

The proofs of (2), (3), (4) are completely elementary and as such 
are left to the reader. Let us just state that for $p\in P(I),
\alpha >\beta$ both in $dom(p)$ we have $p\Vdash$``$\forall i\in 
\omega \setminus n_p\ \dot f_\beta (i)\leq \dot f_\alpha (i)".$ \qed L2
\enddemo

The following Lemma offers an alternative way to view the $P(I),$ 
for $I=\alpha \in Ord.$ It shows how we carefully extend the scale 
of $f_\alpha $'s and it gives an actual computation of the residue 
forcing $P_{\alpha +1}/P(\alpha) .$ We omit the proof, as it is not
 difficult and not very relevant for our purposes.

\proclaim
{Lemma 3} Let $\alpha \in Ord.$ Then $P(\alpha)$ is isomorphic to a dense 
subset of a finite support iteration $\langle R_\beta :\beta \leq \alpha ,
\dot Q_\beta :\beta <\alpha \rangle ,$ where for $\beta \leq \alpha$
\roster
\item $R_\beta \Vdash$``a generic subset of $\dot Q_\beta$ is given by 
a function $\dot f_\beta \in {^\omega \omega }"$
\item $R_\beta \Vdash$``$\delta <\gamma <\beta$ implies that $\dot 
f_\delta \ll \dot f_\gamma "$
\item $R_\beta \Vdash$``$\dot Q_\beta =\{ \langle s,a\rangle : s\in 
{^{<\omega }\omega },a\in [\beta ]^{<\aleph _0}\}$ ordered by 
$\langle s_0,a_0\rangle \geq \langle s_1, a_1\rangle$ if 
$s_0\subset s_1, a_0\subset a_1$ and $\forall i\in 
dom(s_1\setminus s_0)\ \forall \gamma \in a\ s_1(i)\geq 
\dot f_\gamma (i).$ If $G\subset \dot Q_\beta$ is generic 
then $f_\beta$ will be defined as $f_\beta =\bigcup 
\{ s:\langle s,0\rangle \in G \} ."$
\endroster
\endproclaim

The key idea of the proof of the Theorem is given in the following Lemma, 
for which we have not found a satisfactory reference:

\proclaim
{Lemma 4} (Folklore?) Let $\kappa \in Ord .$ Then $\Bbb C(\kappa )\Vdash$
``there is no $\ll$-increasing sequence $\langle \dot f_\alpha :\alpha 
\in (\frak c^+)^V\rangle \subset {^\omega \omega }".$
\endproclaim

\remark
{Remark} This lemma holds true for any antisymmetric relation on 
${^\omega \omega }$ defined from parameters in $V$ in place of $\ll .$
 The only properties of $\Bbb C(\kappa )$ we use are c.c.c. and a form 
of strong homogeneity, see below.
\endremark

\demo
{Proof} Assume for contradiction that we have $q\in C(\kappa ),$ $q\Vdash$``$\langle \dot f_\alpha :\alpha \in (\frak c^+)^V \rangle$ 
is a $\ll$ increasing sequence". By homogeneity of $\Bbb C(\kappa )$ 
we may assume that $q=1.$ Pick $\theta$ large regular and $M_\alpha 
\prec H_\theta ,\alpha \in \frak c^+$ so that $M_\alpha$ are countable
 and $\kappa ,\dot f_\alpha \in M_\alpha .$

\proclaim
{Claim 1} There are $\beta <\alpha <\frak c^+$ and an isomorphism 
$h:\langle M_\alpha ,\in \rangle \to \langle M_\beta ,\in \rangle$ 
such that $h\restriction M_\alpha \cap M_\beta =id$ and 
$h(\dot f_\alpha )=\dot f_\beta .$
\endproclaim

\demo
{Proof of the Claim} By a $\Delta$-system argument, there is 
$X\subset \frak c^+$ of full cardinality and a countable set 
$A$ such that $\langle M_\alpha :\alpha \in X\rangle$ form a 
$\Delta$-system with a root $A.$ Let $A=\{ a_i :i\in \omega \}$ 
be a one-to-one enumeration. There are only $\frak c$ many isomorphism 
types of $\langle M_\alpha ,\in ,\dot f_\alpha ,a_i:i\in \omega \rangle$ 
and so there are $\beta <\alpha \in X$ and $h:\langle M_\alpha ,\in ,
\dot f_\alpha ,a_i:i\in \omega \rangle \to \langle M_\beta ,\in ,\dot
 f_\beta ,a_i:i\in \omega \rangle .$ Obviously $h(\dot f_\alpha )
=\dot f_\beta$ and since $A=\{ a_i:i\in \omega \} =M_\alpha \cap M_\beta$ 
and $h(a_i)=a_i$ we have $h\restriction M_\alpha \cap M_\beta =id$ 
as needed. \qed C1
\enddemo

Let $\beta <\alpha ,h$ be as in the Claim. We extend $h\restriction 
C(\kappa )\cap M_\alpha$ to an automorphism of $C(\kappa )$ and 
$V^{\Bbb C(\kappa )}.$ Let $q\in C(\kappa ), q=q_0\cup q_1 \cup q_2,$ 
where $q_0=q\cap M_\alpha ,$ $q_1=q\cap M_\beta$ and $q_2=q\setminus 
(q_0\cup q_1).$ We set $\bar h(q)=h(q_0)\cup h^{-1}(q_1)\cup q_2 .$ 
Note that $q_0\in M_\alpha$ and $q_1\in M_\beta .$ It can be easily 
verified from the properties of $h$ and the definition of $C(\kappa )$ 
that $\bar h$ is an automorphism of $C(\kappa ).$ As such it can be 
extended to an automorphism of $\Bbb C(\kappa )$ and 
$V^{\Bbb C(\kappa )},$ which by abuse of notation we call $\bar h$ again.

\proclaim
{Claim 2} $\bar h(\dot f_\alpha )=\dot f_\beta ,
\bar h(\dot f_\beta )=\dot f_\alpha .$
\endproclaim

\demo
{Proof of the Claim} Here we use the c.c.c. of $\Bbb C(\kappa ).$ 
It implies that $TC(\dot f_\alpha )\subset M_\alpha .$ One can 
prove by $\in$-induction in $V^{\Bbb C(\kappa )}$ that for any 
$\tau \in M_\alpha ,$ a $\Bbb C(\kappa )$-name such that 
$TC(\tau )\subset M_\alpha ,$ we have $\bar h(\tau )=h(\tau )$
 and $\bar h(h(\tau ))=h^{-1}(h(\tau )).$ \qed C2
\enddemo

The Lemma follows: $\Bbb C(\kappa )\Vdash$``$\dot f_\beta \ll 
\dot f_\alpha "$ and at the same time $\Bbb C(\kappa )\Vdash$
``$\bar h(\dot f_\beta )= \dot f_\alpha \ll \dot f_\beta =\bar 
h(\dot f_\alpha )",$ since $\ll$ is defined by a formula with
 parameters in $V.$ However, this is a contradiction with 
antisymmetricity of $\ll .$ \qed L4
\enddemo

Now we can conclude that $RO(P(\frak c^+))$ witnesses the statement 
of Theorem 1. By Lemma 2 (1) it is semi-Cohen; moreover it adds a 
$\frak c^+$ long $\ll$-increasing sequence of functions in 
${^\omega \omega }$ and so by Lemma 4 it cannot be embedded 
into any $\Bbb C(\kappa ).$ \qed T1
\enddemo

Is it consistent that all semi-Cohen algebras of size $\leq \frak c$ 
can be embedded into a Cohen algebra? Is it implied by, let us say, 
Martin's Axiom?

\subhead
{2. The second example}
\endsubhead

Recently, S. Koppelberg and S. Shelah \cite {KS} proved that subalgebras 
of Cohen algebras are not necessarily Cohen themselves. We have 
simplified their argument and as the forcing used is quite similar 
to the one in Theorem 1, we have found it convenient to include our proof here.

\proclaim
{Theorem 2} \cite {KS} There is a complete subalgebra of 
$\Bbb C(\aleph _2)$ of uniform density $\aleph _2$ which is not 
isomorphic to $\Bbb C(\aleph _2 ).$
\endproclaim

\demo
{Proof} Let $I\subset Ord.$ We define a poset $P(I)$ as follows:

\proclaim
{Definition 4} $P(I)=\{ p:p$ is a function, $dom(p)\in [I]^{<\aleph _0},
\exists n_p \in \omega \ \forall \alpha \in dom(p)$ $p(\alpha )\in {^{n_p}
\omega }\}$ ordered by $p_0\geq p_1$ if $dom(p_0)\subset dom(p_1),$
 $n_{p_0}\leq n_{p_1},$ $\forall \alpha \in dom(p_0)$ $p_0(\alpha )
\subset p_1(\alpha )$ and $\forall i\in n_{p_1}\setminus n_{p_0}\ 
\forall \beta \neq \alpha$ both in $dom(p_0)$ we have $p_1(\beta )
(i)\neq p_1(\alpha )(i).$
\endproclaim

Thus the forcing $P(I)$ is just the standard forcing for adding an 
$I$-indexed family of eventually different functions. By translating 
the proof of Lemma 2 virtually word by word into the new situation, 
we obtain

\proclaim
{Lemma 5}\roster
\item If $J\subset I$ then $P(J)$ is a regular subposet of $P(I);$ 
consequently, $P(I)$ is semi-Cohen.
\item $P(I)$ has uniform density $|I|+\aleph _0.$
\item Let $F\subset P$ be generic. For $\alpha \in I$ we define 
$f_\alpha =\bigcup _{p\in F}p(\alpha )\in {^\omega \omega }.$ 
Then $\alpha \neq \beta \in I$ implies $f_\alpha (i)\neq f_\beta (i)$ 
for all but finitely many $i\in \omega .$
\item $\langle f_\alpha :\alpha \in I\rangle$ determines $F:$ namely, 
$F=\{ p\in P(I):\forall \alpha \in dom(p)\ p(\alpha )\subset f_\alpha$ 
and for all $\alpha >\beta$ both in $dom(p)$ and $i\in \omega 
\setminus n_p$ $f_\beta (i)\leq f_\alpha (i)\} .$ 
\endroster
\endproclaim
 
And as in Theorem 1, we have an iteration representation of $P(I)$ 
for $I=\alpha \in Ord.$ This time we prove it rigorously, as 
it will be instrumental in the proof of Lemma 8.

\proclaim
{Lemma 6} Let $\alpha \in Ord.$ Then $P(\alpha)$ is isomorphic 
to a dense subset of a finite support iteration $\langle R_\beta :
\beta \leq \alpha ,\dot Q_\beta :\beta <\alpha \rangle ,$ where 
for $\beta \leq \alpha$
\roster
\item $R_\beta \Vdash$``a generic subset of $\dot Q_\beta$ is 
given by a function $\dot f_\beta \in {^\omega \omega }"$
\item $R_\beta \Vdash$``$\delta <\gamma <\beta$ implies that 
$\dot f_\delta (i)\neq \dot f_\gamma (i)$ for all but finitely 
many $i\in \omega ".$
\item $R_\beta \Vdash$``$\dot Q_\beta =\{ \langle s,a\rangle : 
s\in {^{<\omega }\omega },a\in [\beta ]^{<\aleph _0}\}$ ordered 
by $\langle s_0,a_0\rangle \geq \langle s_1, a_1\rangle$ if 
$s_0\subset s_1, a_0\subset a_1$ and $\forall i\in dom(s_1\setminus 
s_0)\ \forall \gamma \in a\ s_1(i)\neq \dot f_\gamma (i).$ If 
$G\subset \dot Q_\beta$ is generic then $f_\beta$ will be defined 
as $f_\beta =\bigcup \{ s:\langle s,0\rangle \in G \} ."$
\item $R_\beta \Vdash$``$\dot Q_\beta$ is separative and has 
uniform density $|\beta |+\aleph _0".$
\endroster
\endproclaim

\demo
{Proof} Let $R_\alpha$ be a finite support iteration as described in (3) 
above. We exhibit a dense subset of $R_\alpha$ isomorphic to $P(\alpha) .$ 
Let $D=\{ r\in R_\alpha :\exists n\in \omega \ \forall \beta \in 
dom(r)\ r\restriction \beta \Vdash$``$\dot r(\beta )=\langle 
\check s_{\beta ,r},dom(r)\cap \beta \rangle "$ for some 
$s_{\beta ,r}\in {^n\omega }\} .$ Obviously, the function 
$h:D\to P(\alpha )$ defined by: $h(r)$ is the function $p$ 
for which $dom(p)=dom(r)$ and $p(\beta )=s_{\beta ,r},$ is 
a poset isomorphism.

We must show that $D\subset R_\alpha$ is dense. Fix $r_0\in R_\alpha .$ 
We attempt to construct a sequence $\langle r_i,\beta _i,n_i,s_i,a_i:
i\in \omega \rangle$ so that $r_0\geq r_1\geq \dots$ are elements of 
$R_\alpha ,$ $\beta _0=max(dom(r_0))>\beta _1>\dots$ are ordinals, 
$n_i\in \omega ,s_i\in {^{n_i}\omega },a_i\in [\beta _i]^{<\aleph _0}$ and
\roster
\item"{(a)}" $\forall i\in \omega \ r_{i+1}\restriction \alpha 
\setminus \beta _i=r_i\restriction \alpha \setminus \beta _i$
\item"{(b)}" $\forall i\in \omega \ r_{i+1}\restriction \beta _i
\Vdash$``$ r_i(\alpha _i)=\langle \check s_i,a_i\rangle "$
\item"{(c)}" $\forall i\in \omega \ \beta _{i+1}=max(\beta _i
\cap dom(r_{i+1})).$
\endroster
As the $\beta _i$'s are decreasing, the construction must stop 
at some $m\in \omega .$ Then $\beta _i:i<m$ are a decreasing 
enumeration of $dom(r_m).$ Let $n=max\{ n_i:i<m\} .$ We construct 
$r\in D$ as follows: $dom(r)=dom(r_m)\cup \bigcup \{ a_i:i<m\} ,$ 
$\forall \beta \in dom(r)\ r(\beta )=\langle s_{\beta ,r},dom(r)\cap 
\beta \rangle$ for some $s_{\beta ,r}.$ If $\beta =\beta _i$ for
 some $i<m$ then we will have $s_i\subset s_{\beta ,r}.$ For all 
pairs $(\beta ,k)$ such that $\beta \in dom(r)$ and $k<n$ and either 
$\beta \notin dom(r_m)$ or $\beta =\beta _i$ and $n_i\leq k<n,$ 
we choose the values of $r(\beta )(k)$ to be pairwise distinct and 
also distinct from all elements of $\bigcup _{i<m}rng(s_i).$ By the 
definition of $R_\alpha$ one can easily verify that $r\leq r_m\leq 
r_0.$ As $r\in D,$ we have proven that $D\subset R_\alpha$ is dense. 

(2) and (1) are immediate. We turn to (4) and work in $V^{R_\beta }.$ 
First we prove separativity of $Q_\beta .$ Let $\langle s_0,a_0\rangle
 ,\langle s_1,a_1\rangle \in Q_\beta ,$ $\langle s_0,a_0\rangle \not 
\leq \langle s_1,a_1\rangle .$ We need to find $\langle s_2,a_2\rangle 
\leq \langle s_0,a_0\rangle ,\langle s_2,a_2\rangle \perp \langle s_1,
a_1 \rangle .$ There are three cases:
\roster
\item"{(a)}" $s_1\not \subset s_0.$ Find an extension $\langle s_2,a_0
\rangle \leq \langle s_0,a_0\rangle$ such that there is $i\in dom(s_0)$ 
for which $s_2(i)\neq s_0(i).$ $\langle s_2,a_0\rangle$ will do.
\item"{(b)}" $s_1\subset s_0,a_1\not \subset a_0.$ Fix $\gamma <\beta$ 
such that $\gamma \in a_1\setminus a_0.$ As $f_\gamma$ is eventually 
different from all $f_\delta ,\delta \in a_0,$ we can find $m\in \omega$ 
such that $\forall i>m\ \forall \delta \in a_0\ f_\gamma (i)\neq 
f_\delta (i).$ Consequently we can find $\langle s_2, a_o\rangle \leq 
\langle s_0,a_0\rangle$ so that there is $i\in dom(s_2)\setminus 
dom(s_1)$ with $s_2(i)=f_\gamma (i).$ As $\gamma \in a_1,$ by the 
definition of $Q_\beta$ necessarily $\langle s_2,a_0\rangle \perp 
\langle s_1, a_1\rangle$ and we are done.
\item"{(c)}" $s_1\subset s_0,a_1\subset a_0.$ By the definition of 
$Q_\beta$ in this case either $\langle s_0,a_0\rangle \leq 
\langle s_1,a_1\rangle$ or $\langle s_0,a_0\rangle \perp \langle s_1,a_1
\rangle$ and we are done again.
\endroster

The density part is quite easy. If $|\beta |\leq \aleph _0$ then 
$Q_\beta$ is isomorphic to $\Bbb C(\aleph _0)$ and we are done. 
So assume $|\beta |\geq \aleph _1.$ We have $|Q_\beta |=|\beta |$
 and if $E\subset Q_\beta ,$ $|E|<|\beta |$ then by the pigeonhole 
principle there is $\gamma <\beta$ such that $\forall \langle s,a\rangle
\in E\ \gamma \notin a.$ So $E$ cannot be even locally dense: if 
$\langle s,a\rangle \in Q_\beta$ then $\langle s,a\cup \{ \gamma \} 
\rangle \leq \langle s,a\rangle$ and there is no element of $E$ 
below $\langle s,a\cup \{ \gamma \} \rangle .$ \qed L6
\enddemo

Notice that we have proven a little more than what was required. Namely, 
we have shown that there is a sequence $\langle h_\beta :\beta \leq 
\alpha \rangle$ such that if $\gamma \leq \beta \leq \alpha$ then 
$h_\gamma :RO(P(\gamma ))\to RO(R_\gamma )$ and $h_\gamma :RO(P(\beta ))
\to RO(R_\beta )$ are isomorphisms and $h_\gamma \subset h_\beta ,$ 
if we view $RO(P(\gamma )),RO(R_\gamma )$ canonically as subsets of 
$RO(P(\beta )),RO(R_\beta )$ respectively. The following two lemmas 
finish the proof of Theorem 2. The witness we need is $RO(P(\omega _2)).$

\proclaim
{Lemma 7} $RO(P(\omega _2))$ can be completely embedded into 
$\Bbb C(\aleph _2).$
\endproclaim

\proclaim
{Lemma 8} $RO(P(\omega _2))$ is not isomorphic to $\Bbb C(\aleph _2).$
\endproclaim
 
\demo
{Proof of Lemma 7} Let $R$ be the set of all functions $r$ 
satisfying the following four conditions:
\roster
\item $dom(r)\in [(\omega _2 \times 2)\cup \{ \omega _2\} ]^{<\aleph _0}$
\item $\alpha \in \omega _2, (\alpha ,0)\in dom(r)$ implies 
$r(\alpha ,0)\in {^{<\omega }\omega }$
\item  $\alpha \in \omega _2, (\alpha ,1)\in dom(r)$ implies 
$r(\alpha ,1)\in \omega$
\item $r(\omega _2)$ is a finite function from ${^{<\omega }\omega }$ 
to $\omega$ which is one-to-one on each ${^n\omega },n\in \omega .$
\endroster

$R$ is ordered by coordinatewise extension, i.e. $r_0\geq r_1$ 
if $dom(r_0)\subset dom(r_1)$ and $(\alpha ,0)\in dom(r_0)$ implies 
$r_0(\alpha ,0)\subset r_1(\alpha ,0),$ $(\alpha ,1)\in dom(r_0)$ 
implies $r_0(\alpha ,1)= r_1(\alpha ,1),$ and $r_0(\omega _2)\subset r_1
(\omega _2).$ As the different coordinates behave independently of 
each other, $RO(R)$ is isomorphic to $\Bbb C(\aleph _2).$

\remark
{Motivation} Let $G\subset R$ be generic and work in $V[G].$ For 
$\alpha \in \omega _2$ set $g_\alpha =\bigcup _{r\in G}r(\alpha ,0):
\omega \to \omega ,$ $n_\alpha =r(\alpha ,1)$ for some (any) 
$r\in G,(\alpha ,1)\in dom(r),$ and $h=\bigcup _{r\in G}r(\omega _2):
{^{<\omega }\omega }\to \omega .$ Notice that $g_\alpha$ are mutually 
generic Cohen reals and $h$ is one-to-one on each ${^n\omega },
n\in \omega .$. We need to obtain $F\subset P(\omega _2)$ which is 
$V$-generic. One might be tempted to think that defining $f_\alpha :
\omega \to \omega$ by $f_\alpha (n)=h(g_\alpha \restriction n),$ 
the resulting sequence $\langle f_\alpha :\alpha \in \omega _2\rangle$
 could be $V$-generic for $P(\omega _2).$ (It is indeed a family of 
eventually different functions!) However, this is too simple to work. 
Here the $n_\alpha$'s kick in: if we define $f_\alpha (n)=g_\alpha (n)$ 
for $n< n_\alpha$ and $f_\alpha (n)=h(g_\alpha \restriction n+1)$ for
 $n\geq n_\alpha ,$ the sequence of $f_\alpha$'s will really be 
sufficiently generic. Of course, we have to justify this statement.
\endremark

Let $D$ be the set of those $r\in R$ for which $(\alpha ,0)\in dom(r)$ 
iff $(\alpha ,1)\in dom(r)$ and there is $n_r\in \omega$ such that the 
$r(\alpha ,0)$'s are distinct elements of ${^{n_r}\omega },$ 
$r(\alpha ,1)\leq n_r$ and $dom(r(\omega _2))=\{ \sigma \in 
{^{<\omega }\omega }:\exists \alpha \in \omega _2\ \sigma \subset 
r(\alpha ,0)\} .$ It is not difficult to see that $D\subset R$ is 
dense. We define a function $Proj:D\to P(\omega _2)$ by: $Proj(r)$
 is a function $p\in P(\omega _2)$ such that $\alpha \in dom(p)$ 
iff $(\alpha ,0)\in dom(r),$ $n_p=n_r$ and $p(\alpha )(n)=r(\alpha ,0)(n)$ 
for $n< r(\alpha ,1)$ and $p(\alpha )(n)=r(\omega _2)(r(\alpha ,0)
\restriction n+1)$ for $r(\alpha ,1)\leq n<n_p.$ The following claim 
completes the proof of the Lemma, as it shows that for any $G\subset 
R$ $V$-generic the filter $F=$ upwards closure of $Proj^{\prime \prime }
G\cap D\subset P(\omega _2)$ will be $V$-generic as well.

\proclaim
{Claim 3}\roster
\item For $r_0\geq r_1$ both in $D$ we have $Proj(r_0)\geq 
Proj(r_1)$ in $P(\omega _2).$
\item If $p_0=Proj(r_0)$ and $p_1\leq p_0$ then for some $r_2\leq r_0$ 
in $D$ we have $p_2=Proj(r_2)\leq p_1.$
\endroster
\endproclaim

\demo
{Proof} For (1), let $r_0\geq r_1$ in $D,$ $p_0=Proj(r_0)$ and 
$p_1=Proj(r_1).$ First, by the definition of $Proj,$ $dom(p_0)\subset
 dom(p_1),$ $n_{p_0}\leq n_{p_1}$ and $\forall \alpha \in dom(p_0)\ p_0
(\alpha )\subset p_1(\alpha ).$ The last thing to check here is that for all 
$n$ with $n_{p_0}\leq n<n_{p_1}$ and all $\alpha \neq \beta $ both in 
$dom(p_0)$ we have $p_1(\alpha )(n)\neq p_1(\beta )(n).$ For this note that
\roster
\item"{(a)}" $p_1(\alpha )(n)=r_1(\omega _2)(r_1(\alpha ,0)
\restriction n+1)$ and $p_1(\beta )(n)=r_1(\omega _2)(r_1(\beta ,0)
\restriction n+1)$ (this is because we have $r_1(\alpha ,1)
=r_0(\alpha ,1)<n$ as well as $r_1(\beta ,1)=r_0(\beta ,1)<n$ 
by the definition of $D)$
\item"{(b)}" $r_1(\alpha ,0)\restriction n\neq r_1(\beta ,0)
\restriction n$ (this is since already $r_0(\alpha ,0)
=r_1(\alpha ,0)\restriction n_{r_0}\neq r_1(\beta ,0)\restriction n_{r_0}=r_0(\beta ,0)$ by the definition of $D$ again. Moreover 
$n_{r_0}\leq n)$ and
\item"{(c)}" $r_1(\omega _2)\restriction {^n\omega }$ is one-to-one, 
by (4) of the definition of $R.$
\endroster

For (2), let $p_0=Proj(r_0)$ and $p_1\leq p_0$ in $P(\omega _2).$ 
We have $dom(p_0)\subset dom(p_1)$ and $n_{p_0}\leq n_{p_1}.$ 
We produce $r_2\leq r_0$ in $D$ with $p_2=Proj(r_2)\leq p_1.$ 
We will have $(\alpha ,0)\in dom(r_2)$ iff $\alpha \in dom(p_1)$ 
and $n_{r_2}=n_{p_1}+1.$ 
\roster
\item"{(d)}" For $\alpha \in dom(p_0)$ we choose $r_2(\alpha ,0)$ 
as an arbitrary element of ${^{n_{p_1}+1}\omega}$ extending 
$r_0(\alpha ,0)$ and we set $r_2(\alpha ,1)=r_0(\alpha ,1).$ 
\item"{(e)}" For $\alpha \in dom(p_1\setminus p_0)$ we let 
$r_2(\alpha ,0)\restriction n_{p_1}=p_1(\alpha )$ and choose 
$r_2(\alpha ,0)(n_{p_1})$ to be some integers making all 
$r_2(\alpha ,0), \alpha \in dom(p_1)$ pairwise distinct. 
(Notice that all $r_2(\alpha ,0),\alpha \in dom(p_0)$ as 
defined in (d) are already pairwise distinct, since 
$r_0(\alpha ,0)$ are.) Moreover we set $r_2(\alpha ,1)=n_{p_1}.$ 
\item"{(f)}" Finally, we must define $r_2(\omega _2).$ 
We will have $r_0(\omega _2)\subset r_2(\omega _2)$ and
 for $\sigma \subset r_2(\alpha ,0),$ $\alpha \in dom(p_0),$ 
$n_{p_0}<lth(\sigma )\leq n_{p_1}$ we set $r_2(\sigma )=p_1(\alpha )
(lth(\sigma )-1).$ The reader can easily verify that the requirements 
do not collide and that the resulting function is one-to-one on every 
${^n\omega },n\leq n_{p_1}.$ (For this we use the fact that since 
$p_0\geq p_1$ in $P(\omega _2)$ if $\alpha \neq \beta$ both in 
$dom(p_0)$ and $n_{p_0}\leq n<n_{p_1}$ then $p_1(\alpha )(n)\neq p_1
(\beta )(n).)$ We extend the function we have obtained to $r_2(\omega _2)$
 with $dom(r_2(\omega _2))=\{ \sigma :\exists \alpha \in dom(p_1)\ 
\sigma \subset r_2(\alpha ,0)\}$ and $r_2(\omega _2)$ one-to-one on 
each ${^n\omega },n\leq n_{r_2},$ in an arbitrary manner.
\endroster
Obviously $r_2\in D, r_2\leq r_0.$ Let $p_2=Proj(r_2).$ We need to 
show that $p_2\leq p_1$ in $P(\omega _2).$ By the construction of 
$r_2$ we have that $dom(p_2)=dom(p_1)$ and $\forall \alpha \in 
dom(p_1)\ p_1(alpha )\subset p_2(\alpha ).$ That leaves us with 
the last thing to verify, namely $\forall \alpha \neq \beta$ both 
in $dom(p_1)$ $\forall n\in \omega$ if $n_{p_1}\leq n<n_{p_2}$ 
then $p_2(\alpha )(n)\neq p_2(\beta )(n).$ In our case $n_{p_2}=n_{p_1}+1,$ 
so $n=n_{p_1}\geq r_2(\alpha ,1)$ for any $\alpha \in dom(p_1)$ and 
therefore $p_2(\alpha )(n)=r_2(\omega _2)(r_2(\alpha ,0)).$ 
The required inequality follows from the fact that $r_2(\alpha ,0)$ 
are all distinct (by (e)) and $r_2(\omega _2)\restriction {^{n+1}\omega }$ 
is one-to-one. \qed C3 \qed L7 
\enddemo

\enddemo

\demo
{Proof of Lemma 8} For contradiction, assume that $h:\Bbb C(\aleph _2)
\to RO(P(\omega _2))$ is an isomorphism. By a simple closure argument, 
we can find $\omega _1<\alpha <\omega _2$ such that $h^{\prime \prime}
\Bbb C(\alpha )=RO(P(\alpha) ).$ Fix $G\times H\subset \Bbb C(\alpha )
\times \Bbb C(\omega _2\setminus \alpha )=\Bbb C(\aleph _2)$ generic 
over $V.$ Then $h^{\prime \prime }G\subset RO(P(\alpha) )$ is generic 
over $V$ and is a subset of $h^{\prime \prime }G\times H\subset 
RO(P(\omega _2)),$ which is generic over $V$ again. Note that 
$V[G]=V[h^{\prime \prime }G].$ By mutual genericity of $G$ and $H$ 
we can conclude that $V[G][H]$ is a Cohen extension of $V[G],$ in 
particular every $f:\omega \to \omega$ in $V[G][H]$ comes from an 
extension of $V[G]$ by a $V[G]$-generic filter on $\Bbb C(\aleph _0).$ 
Now $f_\alpha =\dot f_\alpha /h^{\prime \prime }G\times H$ determines 
a $V[G]$-generic filter on the poset $Q_\alpha =\dot Q_\alpha /h^{\prime 
\prime }G\in V[G],$ as computed in Lemma 6. So $V[G]\models$
``$RO(Q_\alpha )$ can be completely embedded into $\Bbb C(\aleph _0)".$ 
However, this is a contradiction as $V[G]\models$``$Q_\alpha$ has 
uniform density $\aleph _1$ by Lemma 6(4) while $\Bbb C(\aleph _0)$ 
has uniform density $\aleph _0."$
\qed L8 \qed T2
\enddemo

\enddemo

\enddocument